\documentclass[11pt]{article}
\usepackage{amssymb,amsmath,latexsym}
\usepackage{color}
\usepackage[mathscr]{eucal}

\newtheorem{thm}{Theorem}[section]

\newtheorem{remark}[thm]{Remark}
\newtheorem{defn}[thm]{Definition}
\newtheorem{cor}[thm]{Corollary}
\newtheorem{example}[thm]{Example}
\newtheorem{lemma}[thm]{Lemma}

\numberwithin{equation}{section}

\def\wt{\widetlde}
\def\pf{\noindent{\bf Proof.} }
\def\qed{{\hfill $\Box$ \bigskip}}

\def\topdot{{\leavevmode
      \raise.2ex\hbox{${\cdot}$}}}

\def\E{{\mathcal E}}
\def\F{{\mathcal F}}

\def\1{{\bf 1}}

\def\bP{{\mbox{\bf P}}}
\def\bE{{\mbox{\bf E}}}
\def\wt{\widetilde}
\def\wbar{\overline}

\def\1{{\bf 1}}
\def\loc{{\rm loc}}

\def\wh{\widehat}
\def\<{\langle}
\def\>{\rangle}

\newcommand{\Rn}{\mathbb{R}^n}

\newcommand{\R}{\mathbb{R}}

\renewcommand{\liminf}{\mathop{\varliminf}}

\providecommand{\MR}{\relax\ifhmode\unskip\space\fi MR }

\providecommand{\href}[2]{#2}

\makeatletter

\@addtoreset{equation}{section} \makeatother

\begin{document}

 \title{\bf On subhamonicity for symmetric Markov processes}

\author{Zhen-Qing Chen\footnote{The research of this author is supported
in part by NSF Grant DMS-0906743.} \quad  and \quad
Kazuhiro Kuwae\footnote{The research of this author is partially supported
 by a Grant-in-Aid for Scientific
Research (C) No.~19540220 from Japan Society for the Promotion of Science }}
\date{(December 12, 2009)}
\maketitle

\begin{abstract}
We establish the equivalence of the analytic and probabilistic notions of subharmonicity in
the framework of general symmetric Hunt processes on locally compact separable metric spaces,
 extending an earlier work
of the first named author on the equivalence of
the analytic and probabilistic notions of harmonicity. As a
corollary, we prove a strong maximum principle for locally bounded
finely continuous
subharmonic functions in the space of functions
locally in the domain of the Dirichlet form under some natural
conditions.
\end{abstract}

\vspace{.6truein}

\noindent {\bf AMS 2000 Mathematics Subject Classification}: Primary 60J45, 31C05;
Secondary  31C25,  60J25.
\bigskip

\noindent {\bf Keywords and phrases:} subharmonic function, uniformly integrable submartingale,
symmetric Hunt process, Dirichlet form, L\'evy system, strong maximum principle

\vfill \eject
\section{Introduction}\label{sec:intro}
It is known that a function being subharmonic in a domain $D\subset \R^d$ can be defined by
$\Delta u\leq0$ on $D$ in the distributional sense; that is, $u\in W^{1,2}_{\loc}(D):=
\{u\in L^2_{\loc}(D)\mid\nabla u\in L^2_{\loc}(D)\}$ so that
$$
\int_{\R^d}\nabla u(x)\cdot\nabla v(x) \,dx\leq0\quad\text{ for any non-negative }v\in C_c^{\infty}(D).
$$
When $u$ is continuous, the above
is equivalent to the following sub-averaging property by running a Brownian motion $X=(\Omega,X_t,\bP_x)_{x\in\R^d}$:
for every relatively compact open subset $U$ of $D$:
$$
u(X_{\tau_U})\in L^1(\bP_x)\quad\text{ and }\quad u(x)\leq\bE_x[u(X_{\tau_U})]\quad\text{ for every }x\in U.
$$
Here $\tau_U:=\inf\{t>0\mid X_t\notin U\}$ is the first exit time from $U$.
A function $u$ is said to be harmonic in $D$ if both $u$ and $-u$ are subharmonic in $D$.
Recently, there
have been interest from several areas of mathematics in determining whether the above two notions
 harmonicity and
subharmonicity remain equivalent in a more general context, such as
symmetric Hunt processes on locally compact separable metric spaces.
For instance, due to their importance in theory and applications, there has been intense interest
recently in studying discontinuous processes and non-local (or integro-differential) operators by both
analytical and probabilistic approaches. See, e.g.~\cite{CK:heatmixed, CK:apriori} and the references therein.
So it is important to identify the connection  between the analytic and probabilistic notions of subharmonic functions. Very recently, in \cite{CHarm} the first named author established the
equivalence between the analytic and probabilistic notions of harmonic functions
for symmetric Markov processes.
Subsequently, the above equivalence is extended in \cite{MZZ}
to non-symmetric Markov processes
associated with sectorial Dirichlet forms.

In this paper, we extend the previous work \cite{CHarm}, that is, we address the question
of the equivalence of the analytic and probabilistic notions of subharmonicity in the context of
symmetric Hunt processes on locally compact separable metric space (Theorem~\ref{T:3.11}).
As a byproduct of our result, we prove that strong maximum principle
  holds for
locally bounded finely continuous
$\E$-subharmonic functions under some conditions (Theorem~\ref{thm:SMP}).
Strong maximum principles for subharmonic functions
 of  second order elliptic operators have been
powerful tools for various fields in analysis and geometry. In \cite{Kw:maximum},
the second named author
established, by using analytic method,
a strong maximum principle for finely continuous $\E$-subharmonic functions
in the framework of irreducible local semi-Dirichlet forms whose Hunt processes satisfy the
absolute continuity condition with respect to the underlying measure, which generalize the classical strong maximum principle for second order elliptic operators (for an extension of
 strong maximum principle
for subharmonicity in the barrier sense, see also \cite{Kw:Calabimax}).
The strong maximum principle developed in \cite{Kw:strongmax, Kw:maximum} can be applied to
 analysis or geometry for geometric singular spaces; Alexandrov spaces or spaces appeared in
 the Gromov-Hausdorff limit of Riemannian manifolds with uniform
 lower Ricci curvature bounds and so on.
 More concretely
in \cite{KwSy:splitting}, we establish splitting theorems for weighted Alexandrov spaces having measure contraction property, which are striking applications of the strong maximum principle treated in
\cite{Kw:strongmax, Kw:maximum} in terms of
symmetric diffusion processes.
The strong maximum principle established in this paper holds for symmetric Markov
processes, which may possibly have  discontinuous sample paths, on locally compact
separable metric spaces,
 should have useful implications
 in the study of  non-local operator or jump type
symmetric Markov processes.

Let $X$ be
be an $m$-symmetric Hunt process on a locally compact separable metric space $E$
whose associated Dirichlet
form $(\E,\F)$ is regular on $L^2(E;m)$. Let $D$
be an open subset of $E$ and $\tau_D$ is the first exit time from $D$
by $X$. Motivated by the example at the beginning of this section, loosely speaking (see next section for
 precise statements), there are two ways to define a function $u$ being subharmonic in $D$ with
 respect to $X$:
(a) (probabilistically) $t\mapsto u(X_{t\wedge \tau_D})$ is a $\bP_x$-uniformly integrable
submartingale for quasi-every $x\in D$; (b)
(analytically) $\E(u,g)\leq0$ for $g\in \F \cap C_c^+(D)$.
We will show in Theorem~\ref{T:3.11} below that these two definitions
are equivalent under some integrability conditions as imposed in the previous work \cite{CHarm}
by the first author.
Note that even in the Brownian motion case, a function $u$ that is subharmonic in $D$
is typically
not in the domain $\F$ of the Dirichlet form. Denote by $\F_{D,\loc}$
the family of functions $u$ on $E$ such that for every
relatively compact open subset $D_1$ of $D$, there is a function $f\in \F$ so that $u=f$
$m$-a.e.~on $D_1$.
To show
these two definitions are equivalent, the crux of the difficulty is to
\begin{enumerate}
\item\label{item:question1} appropriately extend the definition of $\E(u,v)$ to functions $u$ in $\F_{D,\loc}$
that satisfy some minimal integrability
condition when $X$ is discontinuous so that $\E(u,v)$ is well defined for every $v\in \F \cap C_c(D)$;
\item\label{item:question2} show that if $u$ is subharmonic in $D$ in the probabilistic sense, then $u\in \F_{D,\loc}$
and $\E(u,v)\leq0$ for every non-negative $v\in\F \cap C_c(D)$.
\end{enumerate}
The question \eqref{item:question1} is solved in the previous work \cite{CHarm}.
 The main focus of this paper is to address
the second question \eqref{item:question2}. For \eqref{item:question2}, we
 establish a Riesz type decomposition theorem (Lemma~\ref{lem:Riesz}) for $\E$-subharmonic functions,
which is a crucial step in proving our main result.

If one assumes a priori that $u\in \F$, then the equivalence of (a) and (b) is easy to establish.
In next section, we give precise definitions, statements of the main results and their proofs. Four
examples are given to illustrate the main results of this paper.
We
use \lq\lq$:=$\rq\rq as a way of definition. For two real numbers $a$ and $b$, $a\wedge b:=
\min\{a, b\}$.

The results of this paper can be extended to non-symmetric Hunt processes
associated with sectorial Dirichlet forms. We will not pursuit this generalization
here in this paper.

\section{Main result }\label{sec:main}

Let $X=(\Omega, {\mathscr F}_{\infty}, {\mathscr F}_t,
X_t,\zeta,\bP_x, x\in E)$ be an $m$-symmetric right Markov process
on a space $E$,  where  $m$ is a positive $\sigma$-finite
measure  with full topological support on $E$. A cemetery state
$\partial$ is added to $E$ to form $E_\partial:=E\cup\{\partial\}$,
and  $\Omega$ is the totality of right-continuous, left-limited
 sample paths from $[0,\infty)$ to
$E_\partial$ that hold the value $\partial$ once attaining it.
Throughout this paper, every function $f$ on $E$
 is automatically extended to be a function
on $E_{\partial}$
 by setting
 $f(\partial) =0$.
For any $\omega\in\Omega$, we set
$X_t(\omega):=\omega(t)$.  Let $\zeta(\omega):=\inf\{t\geq0\,\mid\,
X_t(\omega)=\partial\}$  be the life time of $X$.
Throughout this paper, we use the convention that $X_\infty (\omega):=\partial$.
 As usual,
${\mathscr F}_{\infty}$ and ${\mathscr F}_t$ are the minimal augmented $\sigma$-algebras
obtained from ${\mathscr F}_{\infty}^0:=\sigma\{X_s\,\mid\, 0\leq s<\infty\}$
and  ${\mathscr F}_t^0:=\sigma\{X_s\,\mid\, 0\leq s\leq t\}$ under $\{\bP_x: x\in E\}$.
For a Borel subset $B$ of $E$, $\tau_B:=\inf\{t\geq0\mid X_t\notin
B\}$ (the {\it exit time} of $B$) is an $(\mathscr F_t)$-stopping time.

The transition semigroup $\{P_t: t\ge 0\}$ of $X$ is defined by
$$
P_tf(x):=\bE_x[f(X_t)]=\bE_x[f(X_t): t< \zeta],\qquad t\ge 0.
$$
Each  $P_t$ may be viewed as an operator on $L^2(E;m)$, and taken as a whole
these operators form a strongly continuous semigroup of
self-adjoint contractions. The Dirichlet form associated with $X$ is the
bilinear form
$$
{\cal E}(u, v):=\lim_{t\downarrow 0}t^{-1}(u-P_tu, v)_m
$$
defined on the space
$$
{\cal F}:=\left\{u\in L^2(E; m)\,\Big|
\,\sup_{t>0}\,\,t^{-1}(u-P_tu, u)_m<\infty \right\}.
$$
Here we use the notation $(f,g)_m:=\int_E f(x)g(x)\, m(dx)$ and we shall use
$|f|_2:=\sqrt{(f,f)_m}$ for $f,g\in L^2(E;m)$.
 $P_t$ is extended to be a strongly continuous semigroup $\{T_t;t\geq0\}$ on $L^2(E;m)$.
Without loss of generality, we may assume that $(\E,\F)$ is a regular Dirichlet form on $L^2(E;m)$ and the $X$ is an $m$-symmetric Hunt process,
where $E$ is a locally compact separable metric space having a one point compactification
$E_{\partial}:=E\cup\{\partial\}$
and $m$ is a positive
Radon measure with full topological support (see \cite{CMR}).

A set $B\subset E_{\partial}$ is called \emph{nearly Borel} if for each probability measure $\mu$ on
$E_{\partial}$, there exist Borel sets $B_1,B_2\subset E_{\partial}$ such that $B_1\subset B\subset B_2$ and $\bP_{\mu}(X_t\in B_2\setminus B_1\text{ for some }t\geq0)=0$. Any hitting time
$\sigma_B:=\inf\{t>0\mid X_t\in B\}$ is an $({\mathscr F}_t)$-stopping time for nearly Borel subset of $E_{\partial}$ (see Theorem~10.7 and the remark after Definition~10.21 in \cite{BG}).
A subset $B$ of $E_{\partial}$ is said to be $X$-invariant if $B$ is nearly Borel and
\begin{align*}
\bP_x(X_t\in B_{\partial}, X_{t-}\in B_{\partial} \ \text{ for all }t\geq0)=1
\qquad \text{ for any }x\in B.
\end{align*}
A set $A$ is \emph{finely open}  if for each $x\in A$ there exists a nearly Borel
subset $B=B(x)$ of $E$ such that $B\supset E\setminus A$ and $\bP_x(\sigma_B>0)=1$.
A set $N$ is called \emph{properly exceptional} if $E\setminus N$ is $X$-invariant and $m(N)=0$.
A nearly Borel set $N$ is called \emph{$m$-polar} if $\bP_m(\sigma_N<\infty)=0$ and any
subset $N$ of $E$ is called \emph{exceptional} if there exists an $m$-polar set $\tilde{N}$ containing $N$.
Clearly any properly exceptional set $N$ is
exceptional.
A function defined q.e.~on an open subset $D$ of $E$ is said to be \emph{q.e.~finely
continuous} on $D$ if there exists a properly exceptional Borel set $N$ such that  $u$ is Borel measurable and finely continuous on $D\setminus N$.
It is known (cf. \cite{FOT}) a quasi-continuous function on $D$ is q.e. finely continuous on $D$.

Let $\F_e$ be the family of $m$-measurable functions $u$ on $E$
such that $\vert u\vert<\infty\,m$-a.e.~and there exists an
${\E}$-Cauchy
sequence $\{u_n\}$ of $\F$ such that $\lim\limits_{n\to\infty}u_n=u\;m$-a.e.
We call $\{u_n\}$ as above an approximating sequence for $u\in\F_e$.
For any $u,v\in\F_e$ and its approximating sequences
$\{u_n\}$, $\{v_n\}$
the limit ${\E}(u,v)=\lim\limits_{n\to\infty}{\E}(u_n,v_n)$
exists and does not
depend on the choices of the approximating sequences for $u$, $v$.
It is known that ${\E}^{1/2}$ on $\F_e$ is a semi-norm and $\F=\F_e\cap L^2(E;m)$.
We call $(\E,\F_e)$ \emph{the extended Dirichlet space} of $(\E,\F)$.
Any $u\in\F_e$ admits a quasi-continuous $m$-version $\tilde{u}$.
Throughout this paper, we always take quasi-continuous $m$-version of the element of $\F_e$, that is,
we omit \emph{tilde} from $\tilde{u}$ for $u\in\F_e$

Let $D$ be an open subset of $E$. We define
\begin{eqnarray*}
\left\{\begin{array}{ll}
\F_D:=\left\{u\in \F\mid
u=0 \ \;\E \hbox{-q.e.~on }E\setminus D\right\}, \\
\E^D(u,v):=\E(u,v) \qquad \hbox{ for }u,v\in\F_D.
\end{array}\right.
\end{eqnarray*}
Then $(\E^D,\F_D)$ is again a regular Dirichlet form on $L^2(D;m)$, which is called the
\emph{part space} in $D$.
Denote by $\F_{D,\loc}$ (resp.~$(\F_D)_{\loc}$) \emph{the space of
functions locally in} $\F$ on $D$ (resp.~\emph{the space of
functions locally in} $\F_D$); that is, $u\in \F_{D,\loc}$
(resp.~$u\in(\F_D)_{\loc}$) if and only if for any relatively
compact open set $U$ with $\wbar{U}\subset D$ there exists
$u_U\in\F$ (resp.~$u_U\in\F_D$) such that $u=u_U$ $m$-a.e.~on $U$.
Note that $(\F_D)_{\loc}\subset \F_{D,\loc}$ and $\1_D\in
(\F_D)_{\loc}$. Any $u\in \F_{D,\loc}$ admits an $m$-version
$\tilde{u}$ of $u$ which is quasi-continuous on $D$. As remarked
above, we always take such $m$-version and omit \emph{tilde} from
$\tilde{u}$ for $u\in \F_{D,\loc}$. We can see that $\F_{D,\loc}\cap
L^{\infty}_{\loc}(D;m)\subset (\F_D)_{\loc}$. Indeed, for $u\in
\F_{D,\loc}\cap L^{\infty}_{\loc}(D;m)$, we can take $u_U\in \F_b$
such that $u=u_U$ $m$-a.e.~on $U$, because
$u_U=(-\|u\|_{U,\infty})\lor u_U\land \|u\|_{U,\infty}$ $m$-a.e.~on
$U$, where $\|u\|_{U,\infty}:=m\text{\rm-ess-}\sup_{U}|u|$. Taking
$\phi\in\F\cap C_c(E)$ with $\phi=1$ on $U$ and $\phi=0$ on $D^c$,
we see $u_U\phi\in\F_D$ and $u=u_U\phi$ $m$-a.e.~on $U$.

\begin{defn}[Sub/Super-harmonicity]
{\rm Let $D$ be an open set in $E$. We say that a nearly Borel
measurable function $u$ defined on $E$ is
\emph{subharmonic (resp.~superharmonic) in $D$} if for any
relatively compact open subset $U$ of $D$ with
  $\wbar{U}\subsetneq D$,
$t\mapsto u(X_{t\wedge\tau_U})$ is a uniformly integrable right
continuous $\bP_x$-submartingale (resp.~$\bP_x$-supermartingale) for
q.e.~$x\in E$. A nearly Borel function $u$ on $E$ is said to be
\emph{harmonic in $D$}
$u$ is  both superharmonic and subharmonic in $D$.
}
\end{defn}

\begin{defn}[Sub/Super-harmonicity in the weak sense]
{\rm Let $D$ be an open set in $E$. We say that a nearly Borel
function $u$ defined on $E$ is \emph{subharmonic (resp.~superharmonic) in $D$
in the weak sense} if $u$ is q.e.~finely continuous in $D$ and
for any relatively compact open subset $U$ with
  $\wbar{U}\subsetneq D$,
$\bE_x[|u|(X_{\tau_U})]<\infty$ for q.e.~$x\in E$ and for q.e.~$x\in
E$, $u(x)\leq\bE_x[u(X_{\tau_U})]$
(resp.~$u(x)\geq\bE_x[u(X_{\tau_U})]$) holds if
$\bP_x(\tau_U<\infty)>0$. A nearly Borel
 measurable
function $u$ on $E$ is said to be \emph{harmonic in $D$ in the weak
sense} if
$u$ is  both superharmonic and subharmonic in $D$ in the weak sense.
}
\end{defn}

Clearly $\1_D$ is superharmonic in $D$ in the weak sense.

\begin{remark}\label{rem:ntegrability}
{\rm
Our definition on the subharmonicity or superharmonicity in the weak sense is different from  what is defined in the Dynkin's textbook \cite{Dyn:Mar} and
 is weaker than it when
$X$ is an $m$-irreducible  diffusion process satisfying
 \eqref{e:2.3}  below.
 Actually,
superharmonicity of $u$ in \cite{Dyn:Mar} requires
 $u$ be locally bounded from below
instead of
the $\bP_x$-integrability of $u(X_{\tau_U})$ for any relatively compact open $U$
with $\wbar{U}\subset D$.
Indeed,  suppose that $X$ is a diffusion process and $u$ is a superharmonic function in $D$
 in the sense of \cite{Dyn:Mar}.
Then for $U$ as above, we have
$$\bE_x[|u (X_{\tau_U})|]\leq \bE_x[u(X_{\tau_U})]+2\bE_x[(-u)^+(X_{\tau_U})]\leq
u(x)+2(-\inf_{\partial U}u)^+<\infty
$$
 for q.e.~$x\in E$.
} \qed
\end{remark}
@
We introduce the following condition:
\begin{equation}
 \text{For any relatively compact open set $U$ with }
 \wbar{U}\subsetneq D,  \
\bP_x(\tau_U<\infty)>0  \  \text{ for q.e.~}x\in U.\label{e:2.3}
\end{equation}
 Condition
\eqref{e:2.3} is satisfied if $(\E,\F)$ is
\emph{$m$-irreducible}, that is, any $(T_t)$-invariant set $B$ is
trivial in the sense that $m(B)=0$ or $m(B^c)=0$.
\smallskip

It will be shown in Lemma \ref{L:3.7} that under condition \eqref{e:2.3},
every subharmonic function in $D$ is a subharmonic function in $D$ in the weak sense.

\medskip

In what follows, all functions denoted by $u$ or $u_i$, $(i=1,2)$
are defined
 on $E$ and are
(nearly) Borel measurable and finite quasi everywhere.

\medskip
For an open set $D\subset E$, we consider the following conditions for a
 (nearly) Borel function $u$ on $E$
 that are introduced in \cite{CHarm}.
For any relatively compact open sets $U,V$ with $\wbar{U}\subset V\subset \wbar{V}\subset D$,
\begin{align}
\int_{U\times (E\setminus V)}|u(y)|J(dxdy)<\infty\label{e:2.1}
\end{align}
and
\begin{align}
\1_U\bE_{\cdot}[(1-\phi_V)|u|(X_{\tau_U})]\in(\F_U)_e,\label{e:2.2}
\end{align}
where $\phi_V\in\F\cap C_c(E)$ with $0\leq\phi_V\leq1$ and $\phi_V=1$ on $V$.

\smallskip
As is noted in \cite{CHarm}, in many concrete cases such as in Examples
2.12-2.14 in \cite{CHarm} (see also Examples \ref{ex:stablelike}-\ref{E:4.4} below),
one can show that condition
\eqref{e:2.1} implies condition \eqref{e:2.2}.

\begin{remark}\label{rem:condition} \rm
\begin{enumerate}
\item
It follows immediately from the proof of \cite[Lemma 2.3]{CHarm} and \cite[Lemma 2.4]{CHarm}
that   condition \eqref{e:2.2} is equivalent to
\begin{equation}
\int_{U\times (E\setminus V) } \bE_x  [(1-\phi_V)|u|(X_{\tau_U})]
 (1-\phi_V(y))|u (y)| J(dx, dy)<\infty.
\end{equation}

\item\label{item:rem1}
In view of  \cite[Lemma~2.3]{CHarm}, every nearly Borel bounded function $u$ on $E$ satisfies
both \eqref{e:2.1} and \eqref{e:2.2}.

\item\label{item:rem2} If $u\in\F_{D,\loc}\cap L^{\infty}_{\loc}(D;m)$, then $u$ is bounded q.e.~on
any relatively compact open $U$ with $\wbar{U}\subset D$, so for any $U,V$ as above, \eqref{e:2.1}
is equivalent to
\begin{align}
\int_{U\times (E\setminus V)}|u(y)-u(x)|J(dxdy)<\infty\label{e:2.1*}
\end{align}
for such $u$.
Clearly, any $u\in\F_e$ satisfies
\begin{align*}
\int_{U\times (E\setminus V)}|u(y)-u(x)|J(dxdy)
\leq J(U\times V^c)^{1/2}\left(\int_{E\times E}|u(y)-u(x)|^2J(dxdy)\right)^{1/2}
<\infty;
\end{align*}
that is, \eqref{e:2.1*} is satisfied by $u\in\F_e$.
Furthermore, by Lemma~2.5 of \cite{CHarm}, both \eqref{e:2.1} and \eqref{e:2.2} hold
for every  $u\in\F_e\cap L^{\infty}_{\loc}(D;m)$. \qed
\end{enumerate}
\end{remark}

The following is
proved in \cite{CHarm}.

\begin{lemma}[cf.~Lemma~2.6 in \cite{CHarm}]\label{L:2.4}
Let $D$ be an open set of $E$. Suppose that
$u$ is a
locally bounded function on $D$ such that
$u$ belongs to $\F_{D,\loc}$ and it satisfies
 condition \eqref{e:2.1}.
Then for every $v\in\F\cap C_c(D)$,
the expression
\begin{align*}
\frac{1}{2}\mu_{\<u,v\>}^c(D)+\frac{1}{2}
 \int_{E\times E}
(u(x)-u(y))(v(x)-v(y))J(dxdy)+\int_Du(x)v(x)\kappa(dx)
\end{align*}
is well-defined and finite; it will still be denoted as $\E(u,v)$.
\end{lemma}

\begin{defn}[$\E$-sub/super-harmonicity]\label{D:2.5}
Let
$u\in\F_{D,\loc}\cap L^{\infty}_{\loc}(D;m)$
 be a function satisfying the condition \eqref{e:2.1}.
We say that $u$ is $\E$-subharmonic (resp.~$\E$-superharmonic) in
$D$ if and only if $\E(u,v)\leq0$ (resp.~$\E(u,v)\geq0$) for every
non-negative $v\in\F\cap C_c(D)$.
A function
$u\in\F_{D,\loc}\cap L^{\infty}_{\loc}(D;m)$
 satisfying
  condition \eqref{e:2.1} is said to be $\E$-harmonic in $D$
if $u$ is both $\E$-superharmonic and $\E$-subharmonic
in $D$. When $D=E$, we omit the phrase \lq in $D$\rq.
\end{defn}

Note that $\1_D\in\F_{D,\loc}$ satisfies \eqref{e:2.1}
 and is $\E$-superharmonic in $D$. It
is $\E$-harmonic in $D$ provided
$\kappa (D)=0$ and $J(D, D^c)=0$.

\medskip
Our main theorem below
is an analogy of Theorem~2.11 in \cite{CHarm} for subharmonic functions.

\begin{thm}\label{T:3.11}
Let $D$ be an open subset of $E$. Suppose that a nearly Borel $u\in L^{\infty}_{\loc}(D;m)$
satisfies conditions \eqref{e:2.1} and \eqref{e:2.2}. Then
\begin{enumerate}
\item\label{item:conclusion1*} $u$ is subharmonic in $D$ if and only if
$u\in (\F_D)_{\loc}$
 and it
is $\E$-subharmonic in $D$.
\item\label{item:conclusion2*} Assume that
 \eqref{e:2.3} holds.
Then $u$ is subharmonic in $D$ if and only if $u$ is subharmonic in
$D$  in the weak sense, that is, for any relatively compact open set
$U$ with $\wbar{U}\subsetneq D$, $u(X_{\tau_U})$ is $\bP_x$-integrable
and $u(x)\leq\bE_x[u(X_{\tau_U})]$ for q.e.~$x\in E$.
\end{enumerate}
\end{thm}

Theorem \ref{T:3.11} will be established through
Lemma \ref{L:3.7} and  Theorems \ref{T:3.8}-\ref{T:3.10}.
As an application of Theorem \ref{T:3.11}, we have the following.

\begin{cor}\label{cor:wedge}
\begin{enumerate}
\item\label{item:wedge1} Let $\eta\in C^1(\R)$ be a convex function
and
 $u\in \F_{D,\loc}\cap L^\infty_{\loc}(D; m)$
be an $\E$-harmonic function in $D$ satisfying
conditions \eqref{e:2.1}--\eqref{e:2.2}.
Suppose that $\eta$ has
 bounded first derivative
or
$u$ is bounded on $E$.
Then  $\eta(u)\in\F_{D,\loc}$ and   is
$\E$-subharmonic
 in $D$ satisfying
conditions \eqref{e:2.1}--\eqref{e:2.2}.

\item\label{item:wedge1b}  The conclusion of \eqref{item:wedge1} remains to true if
$\eta\in C^1(\R)$ is an increasing convex function
and $u\in \F_{D,\loc}\cap L^\infty_{\loc}(D; m)$
is an $\E$-subharmonic function in $D$ satisfying
conditions \eqref{e:2.1}--\eqref{e:2.2}.

\item\label{item:wedge2} Let $p\geq1$ and $u\in \F_{D,\loc}$ be an $\E$-harmonic function in $D$
  that is locally bounded in $D$ and
  satisfies conditions \eqref{e:2.1}--\eqref{e:2.2}. Suppose that $|u|^p$
satisfies conditions \eqref{e:2.1} and \eqref{e:2.2}, and that
\eqref{e:2.3} holds.
Then $|u|^p\in\F_{D,\loc}$ and
   is $\E$-subharmonic in $D$.

\item\label{item:wedge3} Let
  $u_1,u_2\in \F_{D,\loc}\cap L^\infty_\loc (D; m)$
 be $\E$-subharmonic functions in $D$ satisfying conditions  \eqref{e:2.1}--\eqref{e:2.2}.
Then $u_1\lor u_2\in \F_{D,\loc}$ satisfies \eqref{e:2.1}--\eqref{e:2.2} and
 is  $\E$-subharmonic in $D$.
\end{enumerate}
\end{cor}

We say that $X$ satisfies \emph{the absolute
continuity condition with respect to $m$} if the transition
kernel $P_t(x,dy)$ of $X$ is
absolutely continuous with respect to $m(dy)$ for any $t>0$ and
$x\in E$.

As a consequence of Corollary~\ref{cor:wedge}\eqref{item:wedge3}, we have the following strong maximum principle.

\begin{thm}[Strong maximum principle]\label{thm:SMP}
 Assume that $D$ is an open subset of $E$,
$X$ satisfies the absolute continuity condition with respect to $m$ and
$(\E^D,\F_D)$ is $m$-irreducible.
Suppose that $u\in \F_{D,\loc}$ satisfying conditions
\eqref{e:2.1}-\eqref{e:2.2} is a
locally bounded finely
 continuous
$\E$-subharmonic function in $D$. If $u$ attains a maximum at a point $x_0\in D$.
Then $u^+\equiv u^+(x_0)$
 on $D$.
  If in addition
 $\kappa (D)=0$, then $u\equiv u(x_0)$ on $D$.
\end{thm}

 \section{Proofs}

In this section, we present proofs for Theorem~\ref{T:3.11}, Corollary~\ref{cor:wedge} and Theorem~\ref{thm:SMP}. First we prepare a lemma.

\begin{lemma}\label{L:3.1}
For $u\in\F$, the following are equivalent.
\begin{enumerate}
\item $\E(u,v)\leq0$ for every $v\in\F^+$.
\item $T_tu\geq u$ $m$-a.e.  on $E$
 for every $t\geq 0$.
\end{enumerate}
\end{lemma}
\pf
 Clearly (ii) implies (i). The proof of (i)$\Rightarrow$(ii)
is quite similar
 to the proof of Lemma~2.2 in \cite{Kw:maximum}.
 So it is omitted.
Note that we do not assert that
$u\leq 0$ $m$-a.e. on $E$.
\qed

\begin{lemma}\label{L:3.2}
For $u_1,u_2\in\F_e$, if $u_1$ and $u_2$ are $\E$-subharmonic,
 then so is $u_1\lor u_2$.
\end{lemma}

\pf Let $g\in L^1(E;m)$ be such that $0<g\leq1$ $m$-a.e.~on $E$ and  that
$u_1,u_2\in  L^2(E;gm)$.
 The measure $gm$ has full quasi-support with respect to $(\E,\F)$ by
Corollary 4.6.1 in \cite{FOT}.
Denote by $(\wt \E, \wt \F)$ the Dirichlet form of the process $X$ time-changed
by the inverse of $A_t:=\int_0^t g(X_s)ds$. Then by (6.2.22)-(6.2.23) of \cite{FOT},
$(\wt \E, \wt \F_e)=(\E, \F_e)$ and $\wt \F=\wt \F_e\cap L^2(E; gm)$.
By Theorem 6.2.1 of \cite{FOT}, $(\wt \E, \wt \F)$ is a regular Dirichlet
form on $L^2(E; gm)$ with core $\wt \F \cap C_c(E)=\F\cap C_c(E)$.
So $u_1$ and $u_2$ are
$\wt{\E}$-subharmonic functions in $\wt{\F}$.
Let $\{\wt T_t, t\geq 0\}$ be the semigroup associated with $(\wt \E, \wt \F)$.
 From Lemma~\ref{L:3.1}, we see $u_1\leq \wt{T}_tu_1$ and
$u_2\leq\wt{T}_tu_2$ $m$-a.e.~on $E$, which implies $u_1\lor
u_2\leq \wt{T}_t(u_1\lor u_2)$. By Lemma~\ref{L:3.1} again,
$u_1\lor u_2$ is an $\wt{\E}$-subharmonic function in
$\wt{\F}\subset \wt{\F}_e=\F_e$. The conclusion of the lemma now follows. \qed

\begin{lemma}\label{L:3.3}
Let $v_1$  be an excessive function of $X$ and
 $v_2\in \F_e$ such that $v_1\leq v_2$ $m$-a.e. on $E$.
 Then $v_1\in \F_e$ with $\E(v_1, v_1)\leq \E(v_2, v_2)$.
\end{lemma}

\pf As in the proof of Lemma~\ref{L:3.2},
let $g\in
L^1(E;m)$ be such that   $0<g\leq1$ $m$-a.e.~on $E$ and that $v_1,v_2\in
L^2(E;gm)$.  Let $(\wt \E, \wt \F)$ be the time-changed Dirichlet form
with semigroup $\{\wt T_t, t\geq 0\}$ as in the proof of Lemma~\ref{L:3.2}.
Note that $  v_2\in \wt \F_e \cap L^2(E; gm)=\tilde{\F}$.   By
Proposition~2.8 in \cite{BG}, we have $\bE_x[v_1(X_{\tau_t})]\leq
v_1(x)$, where $\tau_t:=\inf\{s>0\mid \int_0^s g(X_u)du>t\}$.
That is, $\wt T_t v_1 \leq v_1$. Observe that  since $\wt  T_t$ is a contraction
operator on $L^2(E; gm)$ for each $t>0$,
$(f, g)\mapsto (f, g-\wt T_t g)_{gm}$ is a non-negative
symmetric quadratic form on $L^2(E; gm)$. Hence
$$
 |(f, \, g-\wt T_t g)_{gm}|\leq  (f, \, f-\wt T_t f)_{gm}^{1/2} \, \cdot \,
  (g, \, g-\wt T_t g)_{gm}^{1/2}.
$$
Since $v_1\in L^2(E; gm)$ and
$$   (v_1,  v_1-\wt T_t v_1)_{gm}\leq   (v_2, v_1-\wt T_t v_1)_{gm}
=(v_1, v_2-\wt T_t v_2)_{gm} \leq  (v_1, v_1-\wt T_t v_1)_{gm}^{1/2 }
\, \cdot \,
  (v_2, v_2-\wt T_t v_2)_{gm}^{1/2 },
$$
we have
$$\lim_{t\to 0} \frac1t  (v_1,\, v_1-\wt T_t v_1)_{gm}\leq
   \lim_{t\to 0} \frac1t  (v_2, \, v_2-\wt T_t v_2)_{gm} = \wt \E (v_2, v_2)<\infty.
$$
It follows that  $v_1\in \wt \F \subset \wt \F_e=\F_e$ with $\E(v_1, v_1)\leq \E (v_2, v_2)$.
  \qed

\begin{lemma}\label{L:3.4}
Let $D$ be an open set of $E$.
 Suppose that $|u_1|\leq |u_2|$ q.e. on $D$ and $u_2$ satisfies   \eqref{e:2.2}.
Then $u_1$ satisfies   \eqref{e:2.2}.
\end{lemma}

\pf  Let $U, V$ be relatively compact open sets such that $\wbar U\subset V \subset \wbar
V\subsetneq D$.
Note that $\bE_x[u(X_{\tau_U})]$ is excessive with respect to
$X^U$ for any non-negative nearly Borel function $u$.
For $i=1, 2$ and
 $v_i(x):=\bE_x[(1-\phi_V)|u_i|(X_{\tau_U})]$, by assumption,
$v_2\in(\F_U)_e$ and $|v_1|\leq |v_2|$ q.e.~on $U$.
 It follows from Lemma~\ref{L:3.3} that
$v_1\in(\F_U)_e$, namely $u_1$ satisfies
\eqref{e:2.2}. \qed

\begin{lemma}[Riesz decomposition]\label{lem:Riesz}
Suppose that $u$ is a non-negative $\E$-superharmonic function in $\F_e$.
Then there exist an $\E$-harmonic function $h\in\F_e$ and a PCAF $A$
  so
 that $u(x)=\bE_x[A_{\zeta}]+h(x)$ q.e.~$x\in E$.  Moreover,
  $t\mapsto u(X_{t })$
 is a uniformly integrable
$\bP_x$-supermartingale for q.e.~$x\in E$.
\end{lemma}

\pf
 There is a
bounded strictly positive $g\in L^1(E; m)$  such that
$u\in L^1(E;gm)\cap L^2(E;gm)$.
 As in the proof of Lemma \ref{L:3.2}, let
 $(\wt \E, \wt \F)$ be time-changed Dirichlet
form of $(\E, \F)$ by the inverse of PCAF $A_t:=\int_0^t g(X_s)ds$.
It is known (cf. \cite{FOT}) that $(\wt \E, \, \wt \F_e)=(\E, \,
\F_e)$ and so $u\in \wt \F_e \cap L^2(E; gm)=\wt \F$. Since
$$ \E(u, \phi) + \int_E u(x) \phi (x) g(x) m(dx)\geq 0
 \qquad \hbox{for every } \phi \in \wt \F^+\cap C_c(E),
$$
by Theorem 2.2.1 of \cite{FOT}, there is a Radon measure $\nu$ so
that
$$ \E(u, \phi) + \int_E u(x) \phi (x) g(x) m(dx)= \int_E \phi(x) \nu
(dx) \qquad \hbox{for every } \phi \in \wt \F\cap C_c(E).
$$
Define $\mu(dx):= \nu (dx) - u(x)g(x) m(dx)$. As $\F=\F_e\cap L^2(E;
m) \subset \wt \F_e \cap L^2(E; gm)=\wt \F$, we have
$$ \E(u,\phi)=\<\mu,\phi\> \qquad \text{ for every }\phi\in \F\cap C_c(E).
$$
Since $u\in \F_e$ is $\E$-subharmonic, the right hand side of the
above display is non-negative. It follows that $\mu$ is a
non-negative Radon measure and consequently it is of finite energy
integral with respect to $(\E, \F)$.
Hence there exists a PCAF $A$ corresponding to $\mu$ such that for
each $\alpha>0$, $u_{\alpha}$ defined by
$u_{\alpha}(x):=\bE_x[\int_0^{\infty}e^{-\alpha t}dA_t]$ is an
element of $\F$  and
\begin{align*}
\E_{\alpha}(u_{\alpha},\phi)=\<\mu,\phi\> \qquad \text{ for every
}\phi\in\F\cap C_c(E),
\end{align*}
(see Theorem 2.2.1 and Lemma 5.1.3 of \cite{FOT}). It is easy to see
that for $0<\alpha<\beta$,
\begin{equation}\label{e:3.1}
 u_\beta -u_\alpha +(\beta -\alpha) R_\alpha u_\beta =0 ,
 \end{equation}
where $\{R_\alpha, \alpha>0\}$ is the resolvent for the process $X$.
 Consequently by \eqref{e:3.1},
for any non-negative $\phi\in L^2(E;m)$,
$F(\beta):=\beta(u_{\beta},\phi)_m$
 satisfies
$ F'(\beta)=(u_{\beta}+R_{\beta}u_{\beta},\phi)_m\geq 0 $.
 This in particular implies that
for $0<\alpha<\beta$, $\alpha u_\alpha \leq \beta u_\beta$ $m$-a.e. on $E$ and so
\begin{align}\label{eq:monotone}
(\alpha u_{\alpha}-\beta u_{\beta} , \, u_{\beta} )_m\leq 0.
\end{align}
 The above then yields
that for $0<\alpha<\beta$
\begin{eqnarray}
0&\leq&\E(u_{\alpha}-u_{\beta},u_{\alpha}-u_{\beta}) \nonumber\\
&=&\E(u_{\alpha},u_{\alpha})-\E(u_{\beta},u_{\beta})+2\E(u_{\beta}-u_{\alpha},u_{\beta})
 \nonumber \\
&=&\E(u_{\alpha},u_{\alpha})-\E(u_{\beta},u_{\beta})+2(
 \alpha u_{\alpha}-\beta u_{\beta}, \, u_{\beta} )_m  \nonumber\\
&\leq& \E(u_{\alpha},u_{\alpha})-\E(u_{\beta},u_{\beta}),
\label{e:3.3}
\end{eqnarray}
which yields the monotone decrease of $\alpha\mapsto
\E(u_{\alpha},u_{\alpha})$.
On the other hand,
\begin{align*}
\E_{\alpha}(u_{\alpha},u_{\alpha}) =\<\mu,\, {u_{\alpha}}\>
=\E(u,u_{\alpha})\leq
\sqrt{\E(u,u)}\sqrt{\E_{\alpha}(u_{\alpha},u_{\alpha})}
\end{align*}
yields
\begin{equation}\label{e:3.4}
\E_\alpha (u_\alpha, \, u_\alpha) \leq \E (u, u) \qquad \hbox{for
every } \alpha >0.
\end{equation}
  Thus the limit $\lim_{\alpha\to0}\E(u_{\alpha},u_{\alpha})$ exists as a finite number.
 Let $\{\alpha_k, k\geq 1\}$ be a decreasing sequence of positive
 numbers that converges to 0. By \eqref{e:3.3}, $\{u_{\alpha_k},
 k\geq 1\}$ is an $\E$-Cauchy sequence in $\F$ and $u_{\alpha_k}$
 converges to $u_0:=\bE_{\cdot}[A_{\zeta}]$ $m$-a.e. on $E$. Hence
 $u_0\in\F_e$.
It follows from \eqref{e:3.4} that
\begin{align*}
 \lim_{\alpha \downarrow 0} |\alpha(u_{\alpha},\phi)_m|
 \leq \lim_{\alpha \downarrow 0}  \alpha \, \| u_{\alpha} \|_{L^2(E; m)} \,
 \|\phi \|_{L^2(E; m)} =0.
\end{align*}
So for every $\phi\in\F\cap C_c(E)$, we get
\begin{align*}
\E(u_0,\phi)&=\lim_{\alpha\to0}\E(u_{\alpha},\phi)
=\lim_{\alpha\to0}\E_{\alpha}(u_{\alpha},\phi)=\<\mu,\phi\>=\E(u,\phi).
\end{align*}
In other words, for $h:=u-u_0\in \F_e$, $\E(h, \phi)=0$ for every
$\phi\in\F\cap C_c(E)$ and hence for every $\phi \in \F_e$. This in
particular implies that $h$ is $\E$-harmonic with $\E(h, h)=0$. By
Lemma~2.2 of \cite{CHarm}, $t\mapsto  h(X_t)$ is a bounded
$\bP_x$-martingale for q.e. $x\in E$.   On the other hand,
\begin{align*}
u_0(X_{t})=\bE_x[A_{\zeta}|\F_{t }]-A_{t }
\end{align*}
is a uniformly integrable $\bP_x$-supermartingale
 for those $x\in E$ such that $u_0(x)=\bE_x[A_{\zeta}]<\infty$.
It follows that
$u(X_t)=
u_0(X_t) +h(X_t)$ is a uniformly integrable
$\bP_x$-supermartingale for q.e.~$x\in E$.
\qed

 \begin{remark}\label{rem:Riesz} \rm
 The assertion of Lemma~\ref{lem:Riesz} also holds
  in the
  quasi-regular Dirichlet form setting.
  In this case, the definition of
 $\E$-superharmonicity of $u\in\F_e$ should be taken to be that
 $\E(u,\phi)\geq0$ for any $\phi\in\F_e^+$. \qed
 \end{remark}

 \begin{lemma}\label{L:3.7}
Let $D$ be an open set and $u$ a nearly Borel function on $E$.
\begin{enumerate}
\item\label{item:equiharm1}
Assume that condition \eqref{e:2.3} holds.  If $u$ is
a subharmonic function
in $D$ and is in $L^2_{\loc}(D;m)$
 then $u$ is subharmonic in $D$ in the weak sense.

\item\label{item:equiharm2} If $u$ is a nearly Borel q.e.~finely
continuous function on $E$ such that  $u$ is subharmonic in $D$ in
the weak sense, then for each relatively compact open set $U$ with
$\wbar{U}\subsetneq D$,
 $\{u(X_{t\wedge \tau_U}) , \, t\geq0\}$
  is a (not necessarily uniformly integrable)  $\bP_x$-submartingale
 for q.e.~$x\in E$.
\end{enumerate}
\end{lemma}

\pf \eqref{item:equiharm1}: Suppose that $u\in L^2_{\loc}(D;m)$ is
subharmonic in $D$.
For any  relatively compact open set  $U$ with
  $\wbar{U}\subsetneq  D$,
by assumption, $\{u(X_{t\wedge\tau_U}),t\geq0\}$ is a uniformly
  integrable $\bP_x$-submartingale for q.e.~$x\in E$.  Then as $t \to
  \infty$,
  $u(X_{t\wedge\tau_U})$
converges in $L^1(\bP_x)$ as well as $\bP_x$-a.s. to some
 random variable
 $\xi$
 for q.e.~$x\in E$. Set
$Y_t:=u(X_{t\wedge\tau_U})$ for $t\in[0,\infty)$ and
$Y_{\infty}:=\xi$.
 Then $\{Y_t,t\in[0,\infty]\}$ is a right-closed $\bP_x$-submartingale for q.e.~$x\in E$.
 Applying the optional sampling theorem (see Theorem~2.59 in \cite{HWY})
 to $\{Y_t,t\in[0,\infty]\}$, we have
 $\bE_x[|u|(X_{\tau_U})]<\infty$  and  $u(x)\leq\bE_x[Y_{\tau_U}]$ for q.e.~$x\in E$.
 Note that $Y_{\tau_U}\1_{\{\tau_U<\infty\}}=u(X_{\tau_U})$ and
$Y_{\tau_U}=u(X_{\tau_U})+\xi \1_{\{\tau_U=\infty\}}$
 $\bP_x$-a.s.~for q.e.~$x\in E$. Set
 $u_2(x):=\bE_x[\xi \1_{\{\tau_U=\infty\}}]$.  We now show that
 $u_2=0$ q.e. on $E$ if
  $\bP_x(\tau_U<\infty)>0$ for  q.e.~$x\in E$.
 It is easy to see that for each $t>0$
$P_t^Uu_2(x)=u_2(x)$ for
 q.e. $x\in U$.
 Note that
\begin{align*}
u_2(x)=\lim_{t\to\infty}\bE_x[u(X_t)\1_U(X_t)\1_{\{\tau_U=\infty\}}]
\end{align*}
 for q.e.~$x\in E$. It follows
 from Schwarz inequality that
 \begin{align*}
 \int_Uu_2(x)^2 m(dx)&\leq \liminf_{n\to\infty}\int_U
  P_n (\1_Uu^2)(x)m(dx)  \leq\int_Uu (x)^2 m(dx)<\infty.
 \end{align*}
Thus $u_2\in\F_U$ and $\E(u_2,u_2)=0$. Applying Lemma~2.2 in
\cite{CHarm} to $u_2$, we have that
 $u_2=0$ q.e. on $U$
if $\bP_x(\tau_U<\infty)>0$ for q.e.~$x\in U$. Therefore we obtain that
$u(x)\leq\bE_x[u(X_{\tau_U})]$
 for q.e. $x\in U$ if
$\bP_x(\tau_U<\infty)>0$ for
q.e.~$x\in U$.
That is, under condition \eqref{e:2.3}, $u$ is subharmonic in $D$ in
the weak sense.

\eqref{item:equiharm2}: Suppose that a nearly Borel q.e.~finely
continuous $u$ is subharmonic in $D$ in the weak sense.
 Then for any
relatively compact open set $U$ with
  $\wbar{U}\subsetneq D$,
$|u (X_{\tau_U})|$ is $\bP_x$-integrable for q.e.~$x\in E$ and for
each $t>0$,
\begin{align} \label{eq:submartingale1}
\bE_x[u(X_{\tau_U})|\F_{t\wedge\tau_U}]\geq u(X_{t\wedge\tau_U}) \qquad \bP_x\hbox{-a.s.}
\end{align}
  for q.e.~$x\in E$.
Set $h_0(x):=\bE_x[u(X_{\tau_U})]$. Then
 $u_0:=h_0-u\geq 0$ q.e. on $U$, $u_0=0$ q.e. on $U^c$, and has the property that
for any relatively compact open subset $O$ with $\wbar{O}\subset U$,
$\bE_x[u_0(X_{\tau_O})]\leq u_0(x)$ for q.e.~on $U$.
By taking a property exceptional set $N$ of $X$ and restricting the
process $X^U$ to $U\setminus N$ if necessary, we have from
Theorem~12.4 in \cite{Dyn:Mar} that the function $u_0$ is excessive
with respect to $X^U$.
 In particular,
 $t\mapsto u_0(X_t)\1_{\{t<\tau_U\}} = u_0(X_{t\wedge \tau_U})$
is a $\bP_x$-supermartingale for
 q.e. $x\in U$.
On the other hand, we see that $\{h_0(X_{t\wedge\tau_U}),t\geq0\}$
is a uniformly integrable $\bP_x$-martingale for
 q.e. $x\in U$.
Therefore
$\{u(X_{t\wedge \tau_U}), \, t\geq 0\}$ is a
$\bP_x$-submartingale
for q.e.~$x\in U$. \qed

The following theorem is an extension of Theorem~2.7 in \cite{CHarm}
to subharmonic functions.

\begin{thm}\label{T:3.8}
Let $D$ be an open set of $E$. Suppose that
  $u\in \F_{D,\loc} \cap L^\infty_\loc (D; m)$
 satisfying conditions
\eqref{e:2.1}-\eqref{e:2.2} is $\E$-subharmonic in $D$.
Then $u$ is subharmonic in $D$.
\end{thm}

\pf
 Let $U$ be a relatively compact open subset of $D$
with $\wbar{U}\subset D$.
Take $\phi\in\F\cap C_c(D)$ such that $0\leq\phi\leq1$ and $\phi=1$ on a relatively compact open neighborhood $V$ of $\wbar{U}$ with $\wbar{V}\subset  D$.
As in the proof of Theorem~2.7 in \cite{CHarm}, we have
  $\phi u\in\F_D$, $h_1(x):=
\bE_x[(\phi u)(X_{\tau_U})]\in\F_e$, $\phi u-h_1\in(\F_U)_e$ and
$$
\E(h_1,v)=0 \qquad  \text{ for every }v\in(\F_U)_e.
$$
 Let
$h_2(x):=\bE_x[((1-\phi)u)(X_{\tau_U})]$, which is well-defined
under condition  \eqref{e:2.2}. Note also that
$\bE_x[|u|(X_{\tau_U})]<\infty$ for q.e.~$x\in E$ under condition
\eqref{e:2.2}. For simplicity, let $h_0:=h_1+h_2$, that
is, $h_0(x):=\bE_x[u(X_{\tau_U})]$.
  By
 the proof of Theorem~2.7 in \cite{CHarm}, we have
$\1_Uh_2\in(\F_U)_e$, $u-h_0\in(\F_U)_e$ and
\begin{align*}
\E(u-h_0,v)\leq0 \qquad \text{ for every   }v\in \F^+\cap C_c(U).
\end{align*}
This in particular implies that $u-h_0$ is $\E$-subharmonic in $U$.
Note that $(u-h_0)^+ \in (\F_U)_e^+ $ and, by Lemma \ref{L:3.2},
$(u-h_0)^+ $ is $\E$-subharmonic in $U$; that is,
\begin{equation}\label{e:3.7}
\E( (u-h_0)^+, \, v)\leq0 \qquad \text{ for every   }v\in \F^+\cap
C_c(U).
\end{equation}
Since $\F\cap C_c(U)$ is $\E$-dense in $(\F_U)_e$, the above
display holds for every non-negative $v\in (\F_U)_e$.
 Indeed, since $(\E, \F_U)$ is a regular Dirichlet form on $L^2(U; m)$,
 for $v\in (\F^+_U)_e$, there is an $\E$-Cauchy sequence $\{v_n, n\geq 1\}$ in
 $\F_U\cap C_c(U)$ that converges to $v$ $m$-a.e. on $U$. By the
 normal contraction property, $\{v_n^+, n\geq
1\}\subset \F^+\cap C_c(U)$ is $\E$-bounded. Thus in view of the
Banach-Saks theorem, there is a subsequence $\{v_{n_k}^+, n\geq 1\}$
whose Ces\'aro mean sequence is $\E$-Cauchy and converges to $v$
$m$-a.e. on $E$. From it we deduce that \eqref{e:3.7} holds for
every $v\in (\F_U)_e^+$. We have in particular
\begin{align}
0\leq \E((u-h_0)^+,(u-h_0)^+)\leq0.\label{eq:vanish}
\end{align}
Thus by Lemma~2.2 in \cite{CHarm}, we get
$(u-h_0)^+(X_t)=(u-h_0)^+(x)$ for all $t\geq0$ $\bP_x$-a.s.~for
q.e.~$x\in E$.
 Consequently,
$(u-h_0)^+(X_t)$ is a bounded $\bP_x$-martingale for q.e.~$x\in E$.
From this fact,  the sets $A:=\{u>h_0\}$ and $A^c=\{u\leq h_0\}$ are
$X$-invariant.
So after taking out a proper exceptional set of $X$ if needed, we
may and do assume that $h$ is finely continuous and that either $A=E$ or
$A^c=E$.

 Suppose $A=E$ and take $x\in A$. Then $u(x)\geq
h_0(x)+\varepsilon$ for some $\varepsilon>0$. We fix such an
$\varepsilon>0$.
 We then have that $u(X_t)\geq h_0(X_t)+\varepsilon$ for all $t\geq0$ $\bP_x$-a.s.
 Consequently,
$$ u(X_{t\wedge\tau_U})\geq
h_0(X_{t\wedge\tau_U})+\varepsilon=\bE_x[u(X_{\tau_U})|\F_{t\wedge\tau_U}]+\varepsilon
\qquad \bP_x \hbox{-a.s.}
$$
 Since
$\bigvee_{t\geq0}\F_{t\wedge\tau_U}=\F_{\tau_U}$ (see (47.7) in
\cite{Shar:gene}), we have
 $u(X_{\tau_U})\geq u(X_{\tau_U})+\varepsilon$
  $\bP_x$-a.s. on $\{\tau_U <\infty\}$ by letting $t\to \infty$.
  This implies that $\bP_x (\tau_U<\infty)=0$ for every $x\in A$.
 Consequently $h_0 = 0$ q.e. on $E$.
  As $u\geq h_0\geq 0$  on $A=E$, we have from above that $u(X_t) =u(X_0)$
  for all $t\geq 0$ $\bP_x$-a.s. for q.e. $x\in E$. This in particular implies
   that  $t\mapsto u(X_{t\wedge \tau_U})$ is a uniformly integrable
  $\bP_x$-martingale for q.e. $x\in E$.

Next suppose $A^c=E$.   Then $h_0-u\in(\F_U)_e$ is a non-negative $\E$-superharmonic function
in $U$.
By Lemma~\ref{lem:Riesz} and Remark~\ref{rem:Riesz},
  $t\mapsto(u-h_0)(X_{t\wedge\tau_U}) $
is a uniformly integrable $\bP_x$-submartingale. By
\eqref{e:2.2},   $\bE_x[|u(X_{\tau_U})|]<\infty$ for
 q.e. $x\in U$,
and so $t\mapsto h_0(X_{t\wedge\tau_U})$ is also a
uniformly integrable $\bP_x$-martingale. This proves that $t\mapsto
u(X_{t\wedge \tau_U})$ is a uniformly integrable
$\bP_x$-martingale. \qed

The following two theorems are the subharmonic counterpart of
Theorem~2.9 in \cite{CHarm}.

\begin{thm}\label{T:3.9}
Let $D$ be an open subset of $E$  and $u$ a nearly Borel
 measurable
function on  $E$ that is locally bounded in $D$.
Suppose one of the following holds:
\begin{enumerate}
\item\label{item:condi1harm} $u$ is subharmonic in $D$.
\item\label{item:condi2harm} $u$ is subharmonic in $D$ in the weak sense
 and \eqref{e:2.3}
    holds.
\end{enumerate}
Then $u\in(\F_{D})_{\loc}$.
\end{thm}

\pf Take a relatively compact open set $U$ with $\wbar{U}\subsetneq
D$. Set $M:=\|u\|_{L^{\infty}(U;m)}$. Then $0\leq M-u\leq2M$ q.e.~on
$U$. If \eqref{item:condi1harm} (resp.~\eqref{item:condi2harm})
holds, then $\{(M-u)(X_{t\wedge\tau_U}),t\geq0\}$
is a uniformly integrable
(resp., by Lemma \ref{L:3.7}(ii), a
   (not necessarily uniformly integrable)) $\bP_x$-supermartingale
for q.e.~$x\in E$.
 Hence for each $t>0$
\begin{align*}
P_t^U(M-u)\leq M-u\ \ \text{ q.e.~on }U.
\end{align*}
 By the same argument as that
after (2.17) in the proof of Theorem~2.9 in \cite{CHarm}, we
 conclude that
$M-u\in \F_{U,\loc}$ and so  $u\in \F_{U,\loc}$.
Since $U$ is arbitrary, we obtain
$u\in\F_{D,\loc}$. Since $u$ is locally bounded on $D$,
this implies that  $u\in(\F_{D})_{\loc}$. \qed

\begin{thm}\label{T:3.10}
Let $D$ be an open subset of $E$ and $u$ be a nearly Borel function on $E$
 that is in $\F_{D, \loc}\cap L^\infty_\loc (D; m)$ and satisfies
 conditions
\eqref{e:2.1} and \eqref{e:2.2}. Suppose one of the
following holds:
\begin{enumerate}
\item\label{item:condi1harm*} $u$ is subharmonic in $D$.
\item\label{item:condi2harm*} $u$ is subharmonic in $D$ in the weak sense and
\eqref{e:2.3} holds.
\end{enumerate}
Then $u$ is $\E$-subharmonic in $D$.
\end{thm}

\pf  Note that $u$ is automatically q.e.~finely continuous
 in $D$.
In either case,
 by the assumption and Lemma \ref{L:3.7}(ii),
for any relatively compact open set $U$ with $\wbar{U}\subsetneq D$,
we have $\bE_x[|u (X_{\tau_U})|]<\infty$ for q.e.~$x\in E$.
 Take $\phi\in\F\cap C_c(D)$ with $0\leq\phi\leq1$ and $\phi=1$
on a relatively compact open set $V$ with $\wbar{U}\subset V\subset
\wbar{V}\subsetneq D$. Set $h_1(x):=\bE_x[(\phi u)(X_{\tau_U})]$ and
$h_2(x):= \bE_x[((1-\phi) u)(X_{\tau_U})]$, which is
q.e.~well-defined as $\bE_x[|u|(X_{\tau_U})]<\infty$ for q.e.~$x\in
E$.
 By the same argument as that for
Theorems~2.9 and 2.7 in \cite{CHarm}, we see that $\phi u\in\F_D$,
$h_1\in(\F_D)_e$, $\1_Uh_2\in(\F_U)_e$, $h_2=\1_Uh_2+u-\phi u\in
\F_{U,\loc}$ and $\E(h_1,v)=\E(h_2,v)=0$ for any $v\in(\F_U)_e$.
Therefore $h_0(x):=h_1(x)+h_2(x)=\bE_x[u(X_{\tau_U})]$ satisfies
$u_0:=h_0-u=\1_Uh_2+h_1-\phi u\in(\F_U)_e$. For the case
\eqref{item:condi2harm*}, as in the proof of Lemma~\ref{L:3.7} we
see $u_0$ is excessive with respect to the subprocess $X^U$. For the case
\eqref{item:condi1harm*}, we have the same conclusion easily. Then
for each $n\in\mathbb{N}$, we have
\begin{align*}
P_t^U(u_0\wedge n)(x)\leq (u_0\wedge n)(x) \qquad \text{ for
q.e.}~x\in U.
\end{align*}
 Since $u_0\wedge n\in \F_U$ because $m(U)<\infty$,
Lemma~\ref{L:3.1} leads us to
\begin{align*}
\E(u_0\wedge n,\phi)\geq0 \qquad \text{ for every   }\phi\in\F^+\cap
C_c(U).
\end{align*}
On the other  hand, $\{u_0\wedge n\}$ is an $\E$-bounded sequence.
There is a subsequence of $\{u_0\wedge n\}$ whose Ces\'aro mean
sequence is $\E$-Cauchy, and so is $\E$-convergent to $u_0$. We thus
have $\E(u_0, \phi)\leq 0$ for every $\phi \in \F^+\cap C_c(U)$, and
so
\begin{align*}
\E(u,\phi)\leq\E(h_0,\phi)=\E(h_1+h_2,\phi)=0 \qquad \text{ for
every   }\phi\in\F^+\cap C_c(U).
\end{align*}
Since $U$ is arbitrary, we obtain the $\E$-subharmonicity of $u$ in $D$.
\qed

\noindent{\bf Proof of Theorem~\ref{T:3.11}.}
Theorem~\ref{T:3.11} is an easy consequence
 of Lemma \ref{L:3.7},
Theorems~\ref{T:3.8}, \ref{T:3.9} and \ref{T:3.10}.  \qed

\noindent{\bf Proof of Corollary~\ref{cor:wedge}.}
\eqref{item:wedge1}: By Theorem~\ref{T:3.8}, for each relatively
compact open set $U$ with $\wbar{U}\subsetneq D$,
$\{u(X_{t\wedge\tau_U}),t\geq0\}$ is a uniformly integrable
$\bP_x$-martingale for q.e.~$x\in E$. First  assume that $\eta$
has bounded first derivative. Since
 $|\eta(t)-\eta(s)|\leq \sup_{\ell\in\R}|\eta'(\ell)|\cdot|t-s|$ for $t,s\in\R$,
 $\eta(u)\in\F_{D,\loc}$. Meanwhile, $|\eta(u)|\leq \sup_{\ell\in\R}|\eta'(\ell)||u|+|\eta(0)|$
yields that  $\{\eta(u)(X_{t\wedge\tau_U}),t\geq0\}$ is uniformly
integrable under $\bP_x$ for q.e.~$x\in U$   and
$\eta(u)$ satisfies \eqref{e:2.1}--\eqref{e:2.2} by
Lemma~\ref{L:3.4}. (Recall that any bounded function satisfies
\eqref{e:2.1}--\eqref{e:2.2}.) By Jensen's inequality
$\{\eta(u)(X_{t\wedge\tau_U}),t\geq0\}$ is a $\bP_x$-submartingale for q.e. $x\in U$.
 The  $\E$-subharmonicity of $\eta(u)$ in $D$  now follows from
Theorem~\ref{T:3.10}. Next we assume
the boundedness of $u$ on $E$.
Then $\eta(u)\in \F_{D,\loc}$ is bounded on $E$
and it satisfies
\eqref{e:2.1}--\eqref{e:2.2}. The rest  of the proof
is similar as above.

\par  \eqref{item:wedge1b}: The proof is the same as that for  \eqref{item:wedge1}.

\par  \eqref{item:wedge2}:
By Theorem \ref{T:3.11},
$\bE_x[|u(X_{\tau_U})|]<\infty$ and $u(x)=\bE_x[u(X_{\tau_U})]$ for
q.e.~$x\in E$,
 and consequently
$u(X_{t\wedge\tau_U})=\bE_x[u(X_{\tau_U})|\F_{t\wedge\tau_U}]$ for
q.e.~$x\in E$.
 Since $u\in L^\infty_\loc (D; m)$,  $|u|^p\in\F_{D,\loc}\cap L^\infty_\loc (D; m)$.
 Therefore for every  $\phi\in\F\cap C_c(D)$ with
$0\leq\phi\leq1$ and
 $\phi=1$ on an open neighborhood $V$ of $\wbar{U}$ with $\wbar{V}\subsetneq D$,
 $\bE_x[\phi|u|^p(X_{\tau_U})]<\infty$ for q.e.~$x\in E$.
  By assumption,  $\bE_x[(1-\phi)|u|^p(X_{\tau_U})]<\infty$ for q.e.~$x\in E$.
 Therefore $\bE_x[|u|^p(X_{\tau_U})]<\infty$ for q.e.~$x\in E$.
 By  Jensen's inequality, $|u|^p$ is subharmonic in $D$ in the weak sense.
The $\E$-subharmonicity of $|u|^p$ in $D$ now follows from
 Theorem~\ref{T:3.10}.

 \par
 \eqref{item:wedge3}: Note that $|u_1\vee u_2|\leq|u_1|+|u_2|$.
 So by  Lemma~\ref{L:3.4},
 $u_1\vee u_2$ satisfies conditions \eqref{e:2.1}--\eqref{e:2.2}.
 The conclusion follows from Theorem \ref{T:3.11}.
\qed

\noindent{\bf Proof of Theorem~\ref{thm:SMP}.} Since
$u^+(x_0)\geq0$ and $\1_D\in\F_{D,\loc}$ is $\E$-superharmonic in
$D$, $u^+(x_0)-u\in \F_{D,\loc}$ is a finely continuous Borel
measurable non-negative $\E$-superharmonic function in $D$.  Hence
so is $v:=u^+(x_0)-u^+=(u^+(x_0)-u)\wedge u^+(x_0)$ by
Corollary~\ref{cor:wedge}\eqref{item:wedge3}. We set $Y:=\{x\in
D\mid v(x)>0\}$. By Theorem~\ref{T:3.8}, $v$ is also excessive with
respect to $X^D$, so is $\1_Y$ (cf.~\cite{Kw:strongmax}). In
particular, $\1_Y$ is finely continuous with respect to $X^D$. By
Theorem~5.3 in \cite{Kw:maximum}, the irreducibility of $(\E^D,\F_D)$
implies the connectedness of the fine topology on $D$ induced by the part process
 $X^D$.
Thus either $Y=\emptyset$ or $D\setminus Y=\emptyset$. Since $x_0\in
D\setminus Y$, we have $Y=\emptyset$.
So $u^+\equiv u^+(x_0)$ on $D$.
The proof for the case
$\kappa(D)=0$
is quite similar, so it is omitted.
\qed

\section{Examples}

\begin{example}[Stable-like process on $\R^d$]\label{ex:stablelike}
{\rm Consider the following Dirichlet form $(\E,\F)$ on $L^2(\R^d)$, where
\begin{align*}
\left\{\begin{array}{cc}\F=W^{\alpha/2,
}(\R^d)=\left\{u\in L^2(\R^d)\,\Bigl|\,
\int_{\R^d\times\R^d}{(u(x)-u(y))^2}{|x-y|^{d+\alpha}}dxdy<\infty \right\},  \\ \E(u,v)=\frac{1}{2}\int_{\R^d\times\R^d}{(u(x)-u(y))(v(x)-v(y)}{|x-y|^{d+\alpha}}c(x,y)dxdy\text{ for }u,v\in\F.\end{array}\right.
\end{align*}
Here $d\geq1$, $\alpha\in]0,2[$, and $c(x, y)$ is a symmetric function in $(x,y)$ that is bounded between two positive
constants. In literature, $W^{\alpha/2,2}(\R^d)$ is called the Sobolev space on $\R^d$ of fractional order $(\alpha/2, 2)$. For an open
set $D\subset \R^d$, $W^{\alpha/2,2}(D)$ is similarly defined as above but with $D$ in place of
$\R^d$. It is easy to check that $(\E,\F)$
is a regular Dirichlet form on $L^2(\R^d)$ and its associated symmetric Hunt process $X$ is called symmetric
$\alpha$-stable-like process on $\R^d$, which is studied in \cite{CK:heatmixed}.
When $c(x, y)\equiv A(d,-\alpha):=\frac{\alpha 2^{d+\alpha}\Gamma(\frac{d+\alpha}{2})}{2^{d+1}\pi^{d/2}\Gamma(1-\frac{\alpha}{2})}$, the process $X$ is nothing but
the rotationally symmetric $\alpha$-stable process on $\R^d$.
It is shown in \cite{CK:heatmixed} that the symmetric $\alpha$-stable-like process
$X$ has strictly positive jointly continuous
transition density function $p_t(x, y)$
with respect to the Lebesgue measure on $\R^d$
and hence is irreducible. Moreover, there is constant $c > 0$ such that
\begin{align}
p_t(x,y)\leq ct^{-d/\alpha} \ \ \ \text{ for }t>0\text{ and } x,y\in\R^d .
\end{align}
 Consequently, by \cite[Theorem]{ChungGreen},
\begin{align}
\sup_{x\in U}\bE_x[\tau_U]<\infty.\label{eq:greenbdd}
\end{align}
for any open set $U$ having finite Lebesgue measure.
Note that in this example, the jumping measure
$$
J(dxdy) = \frac{c(x, y)}{
|x-y|^{d+\alpha}}dxdy
$$
Hence for any non-empty open set $D\subset \R^d$, condition \eqref{e:2.1}
is satisfied if and only if $(1\wedge |x|^{-d-\alpha})u(x)
\in L^1(\R^d)$
(or equivalently, $u(x)/(1+|x|)^{d+\alpha}\in L^1(\R^d)$).
 As is shown in \cite[Example 2.12]{CHarm}, condition \eqref{e:2.2} is automatically
 satisfied for such $u$.
 When $\alpha \in ]1,2[$, every (globally) Lipschitz function $u$ on $\R^d$
satisfies the condition \eqref{e:2.1}, that is, $(1\wedge |x|^{-d-\alpha})u(x)\in L^1(\R^d)$ holds.
Consequently \eqref{e:2.2} holds for any Lipschitz function $u$
provided $\alpha\in]1,2[$.
Indeed, for any relatively compact open sets $U$, $V$ with $\wbar{U}\subset V\subset\wbar{V}\subset D$,
\begin{align*}
\int_{U\times V^c}\frac{|u(y)-u(x)|}{|x-y|^{d+\alpha}}dxdy
&\leq \|u\|_{\rm\scriptsize Lip}
\int_{U\times V^c}\frac{|x-y|}{|x-y|^{d+\alpha}}dxdy \\
&\leq \|u\|_{\rm\scriptsize Lip}
\sigma(\mathbb{S}^{d-1})\int_U\int_{d(x,V^c)}^{\infty}r^{-\alpha}drdx\\
&\leq\|u\|_{\rm\scriptsize Lip}
|U|\sigma(\mathbb{S}^{d-1})\frac{d(U,V^c)^{1-\alpha}}{\alpha-1}<\infty,
\end{align*}
 and so by Remark~\ref{rem:ntegrability}, \eqref{e:2.1} holds.
Here $\|u\|_{\rm\scriptsize Lip}:=\sup_{x,y\in \R^d}\frac{|u(x)-u(y)|}{|x-y|}$, $|U|$ denotes the volume of $U$ and $\sigma(\mathbb{S}^{d-1})$ is the $(d-1)$-dimensional volume of the unit sphere $\mathbb{S}^{d-1}$.

Theorem~\ref{T:3.11} says that for an open set $D$ and a nearly Borel
 function $u$ on $\R^d$ that is locally bounded on $D$
with $(1\wedge |x|^{-d-\alpha})u(x)
\in L^1(\R^d)$, the following are equivalent.
\begin{enumerate}
\item\label{item:subharmstable1} $u$ is subharmonic in $D$;
\item\label{item:subharmstable2} For every relatively compact open subset $U$ of $D$, $u(X_{\tau_U})\in L^1(\bP_x)$ and
$u(x)\leq\bE_x[u(X_{\tau_U})]$ for q.e.~$x\in U$;
\item\label{item:subharmstable3} $u\in \F_{D,\loc} =
W^{\alpha/2 ,2}_{\loc}(D)$ and
\end{enumerate}
$$
\int_{\R^d\times\R^d}
(u(x)-u(y))(v(x)-v(y))\frac{c(x,y)}{|x-y|^{d+\alpha}}dxdy\leq0\text{ for every }v\in
W^{\alpha/2,2}(D)\cap C_c^+(D).
$$

}
\end{example}

\begin{example}[Symmetric Relativistic $\alpha$-stable Process]\label{ex:relastable}
\rm Take $\alpha\in]0,2[$ and $m\geq0$.
Let $\text{X}^{\text{\tiny\rm R},\alpha}=(\Omega,X_t,\bP_x)_{x\in\R^d}$ be a
L\'evy process on $\R^d$ with
  \[ \bE_{0}\bigl[e^{i \<  \xi, X_{t}  \>}\bigr]=
e^{-t((|\xi|^2+m^{2/\alpha })^{\alpha/2}-m)}.\]
If $m>0$, it is called the \emph{relativistic $\alpha$-stable process with mass $m$} (see \cite{Ryz:relativisticstable}).
In particular, if $\alpha=1$ and $m>0$, it is called the \emph{relativistic free Hamiltonian process} (see \cite{HS:PP}).
When $m=0$, $\text{X}^{\text{\tiny\rm R},\alpha}$ is nothing but the usual \emph{symmetric
$\alpha$-stable process}.
Let $(\E^{\text{\tiny\rm R},\alpha},\F^{\text{\tiny\rm R},\alpha})$
be the Dirichlet form on $L^2(\R^d)$ associated with $\text{\bf X}^{\text{\tiny\rm R},\alpha}$.
Using Fourier transform $\wh{f}(x):=\frac{1}{(2\pi)^{d/2}}\int_{\R^d}e^{i\<x,y\>}f(y)dy$, it follows
from Example~1.4.1 of \cite{FOT} that
\begin{align*}
\left\{\begin{array}{cc}\F^{\text{\tiny\rm R},\alpha}&:=\displaystyle{\left\{f\in L^2(\R^d)\;\Bigl|\; \int_{\R^d}
|\wh{f}(\xi)|^2\left((|\xi|^2+m^{2/\alpha})^{\alpha/2} -m \right)d\xi<\infty  \right\}}, \\
\E^{\text{\tiny\rm R},\alpha}(f,g)&:=\displaystyle{\int_{\R^d}
\wh{f}(\xi)\bar{\wh{g}}(\xi)\left((|\xi|^2+m^{2/\alpha})^{\alpha/2} -m \right)d\xi \quad\text{ for }f,g\in \F^{\text{\tiny\rm R},\alpha}}.
\end{array}\right.
\end{align*}
It is shown by Ryznar \cite{Ryz:relativisticstable} that
the semigroup kernel $p_t(x,y)$ of $\text{\bf X}^{\text{\tiny\rm R},\alpha}$ is given by
\[\hspace{0cm}p_t(x,y)=e^{mt}
\int_0^{\infty}\left(\frac{1}{4\pi s}\right)^{d/2}e^{-\frac{|x-y|^2}{4s}}e^{-sm^{\frac{2}{\alpha}}}\theta_{\frac{\alpha}{2}}(t,s)ds
,\]
where $\theta_{\delta}(t,s)$ is the nonnegative function called the \emph{subordinator}
whose Laplace transform is given by
$$
\int_0^{\infty}e^{-\lambda s}\theta_{\delta}(t,s)ds=e^{-t\lambda^{\delta}}.
$$
Then we see the conservativeness of $\text{\bf X}^{\text{\tiny\rm R},\alpha}$ and the irreducibility of
$(\E^{\text{\tiny\rm R},\alpha},\F^{\text{\tiny\rm R},\alpha})$.
From Lemma~3 in \cite{Ryz:relativisticstable}, there exists $C(d,m)>0$ depending only on $m$ and $d$  such that
$$
\sup_{x,y\in\R^d}p_t(x,y)\leq C(d,m)\left(m^{d/\alpha-d/2}t^{-d/2}+t^{-d/\alpha}\right)\quad\text{ for any }t>0.
$$
This yields by \cite[Theorem 1]{ChungGreen} that \eqref{eq:greenbdd} holds
for any open set $U$ having finite Lebesgue measure.
It is shown in \cite{CS3} that the corresponding jumping measure
satisfies
\begin{align*}
J(dxdy)=\frac{c(x,y)}{|x-y|^{d+\alpha}}dxdy\ \ \text{ with }\ \ c(x,y):=\frac{A(d,-\alpha)}{2}
\Psi(m^{1/\alpha}|x-y|),
\end{align*}
where $A(d,-\alpha)=\frac{\alpha 2^{d+\alpha}\Gamma(\frac{d+\alpha}{2})}{2^{d+1}\pi^{d/2}\Gamma(1-\frac{\alpha}{2})}$, and the function $\Psi$ on $[0,\infty[$ is given by
$\Psi(r):=I(r)/I(0)$ with $I(r):=\int_0^{\infty}s^{\frac{d+\alpha}{2}-1}e^{-\frac{s}{4}-\frac{r^2}{s}}ds$.
 Note that
$\Psi$ is decreasing and satisfies
$\Psi(r)\asymp e^{-r}(1+r^{(d+\alpha-1)/2})$ near $r=\infty$, and $\Psi(r)=1+\Psi''(0)r^2/2+o(r^4)$ near $r=0$.  In particular, for $b>0$ we have
\begin{align}
0<\inf_{r\geq0}\frac{\Psi(m^{1/\alpha}(r+b))}{\Psi(m^{1/\alpha}r)}\leq\sup_{r\geq0}
\frac{\Psi(m^{1/\alpha}(r+b))}{\Psi(m^{1/\alpha}r)}<\infty\label{eq:asymprelstable}
\end{align}
and
\begin{align*}
\left\{\begin{array}{cc}\F^{\text{\tiny\rm R},\alpha}&\hspace{-2cm}=\displaystyle{\left\{f\in L^2(\R^d)\;\Bigl|\;\int_{\R^d\times\R^d}|f(x)-f(y)|^2\frac{c(x,y)}{|x-y|^{d+\alpha}}dxdy<\infty  \right\}}, \\
\E^{\text{\tiny\rm R},\alpha}(f,g)&=\displaystyle{\int_{\R^d\times\R^d}
(f(x)-f(y))(g(x)-g(y))\frac{c(x,y)}{|x-y|^{d+\alpha}}dxdy
\quad\text{ for }f,g\in \F^{\text{\tiny\rm R},\alpha}}.
\end{array}\right.
\end{align*}
Applying \eqref{eq:asymprelstable}, we can obtain that
for any relatively compact open sets $U, V$ with $0\in U$ and
$\wbar{U}\subset V\subset \wbar{V}\subset D$, condition \eqref{e:2.1} is satisfied if and only if
$\Psi(m^{1/\alpha}|x|)(1\wedge |x|^{-d-\alpha})u(x)\in L^1(\R^d)$ (equivalently $\Psi(m^{1/\alpha}|x|)u(x)/(1+|x|)^{d+\alpha}\in L^1(\R^d)$).
Similarly,
any function
$u$ with $\Psi(m^{1/\alpha}|x|)(1\wedge |x|^{-d-\alpha})u(x)\in L^1(\R^d)$ also
satisfies the condition \eqref{e:2.2}
in the same way as in Example~\ref{ex:stablelike}.
For $u\in L^{\infty}_{\loc}(D;m)\cap \F_{D,\loc}^{\text{\tiny\rm R},\alpha}$,
we can deduce  \eqref{e:2.1} and \eqref{e:2.2} if
$\Psi(m^{1/\alpha}|x|)\left(1\wedge|x|^{-d-\alpha}\right)u(x)\in L^1(\R^d)$ without assuming
$0\in U$.
 In this case,
 \eqref{e:2.1} for any relatively compact open $U,V$ with $\wbar{U}\subset V\subset \wbar{V}\subset D$ is equivalent to $\Psi(m^{1/\alpha}|x|)\left(1\wedge|x|^{-d-\alpha}\right)u(x)\in L^1(\R^d)$.
Moreover, any (globally) Lipschitz function $u$ satisfies
\eqref{e:2.1}, consequently \eqref{e:2.2} holds for such $u$.
Indeed, for any relatively compact open sets $U$, $V$ with $\wbar{U}\subset V$,
\begin{align*}
\int_{U\times V^c}\frac{|u(y)-u(x)|}{|x-y|^{d+\alpha}}c(x,y)dxdy
&\leq \frac{A(d,-\alpha)}{2}\|u\|_{\rm\scriptsize Lip}
\int_{U\times V^c}\frac{|x-y|\Psi(m^{1/\alpha}|x-y|)}{|x-y|^{d+\alpha}}dxdy \\
&\leq\frac{A(d,-\alpha)}{2} \|u\|_{\rm\scriptsize Lip}
\sigma(\mathbb{S}^{d-1})\int_U\int_{d(x,V^c)}^{\infty}\Psi(m^{1/\alpha}r)r^{-\alpha}drdx\\
&\leq
 C
\int_{d(U,V^c)}^{\infty}\!\!\!e^{-m^{1/\alpha}r}(1+m^{\frac{d+\alpha-1}{2\alpha}}
r^{\frac{d+\alpha-1}{2}})r^{-\alpha}dr<\infty,
\end{align*}
 and so
\eqref{e:2.1} holds by Remark~\ref{rem:ntegrability}.
Here $C$ is a positive constant.

 By Theorem \ref{T:3.11},
for an open set $D$ and a nearly Borel
 function $u$ on $\R^d$ that is locally bounded on $D$
 with $\Psi(m^{1/\alpha}|x|)(1\wedge |x|^{-d-\alpha})u(x)
\in L^1(\R^d)$, the following are equivalent.
\begin{enumerate}
\item\label{item:subharmstable1b} $u$ is subharmonic in $D$;
\item\label{item:subharmstable2b}
For  every relatively compact open subset $U$ of $D$, $u(X_{\tau_U})\in L^1(\bP_x)$ and
$u(x)\leq\bE_x[u(X_{\tau_U})]$ for q.e.~$x\in U$;
\item\label{item:subharmstable3b}
$u\in \F^{\text{\tiny\rm R},\alpha}_{D,\loc}$ and
\end{enumerate}
$$
\int_{\R^d\times\R^d}
(u(x)-u(y))(v(x)-v(y))\frac{\Psi(m^{1/\alpha}|x-y|)}{|x-y|^{d+\alpha}}dxdy\leq0
\qquad \text{for every } v\in
\F^{\text{\tiny\rm R},\alpha}_{D}\cap C_c^+(D).
$$

\medskip

One may ask concrete examples of $\E$-(sub/super)-harmonicity on $D$.
To answer this question, in what follows, we assume $d>2$ ($d>\alpha$ if $m=0$).
Applying Theorems~3.1 and 3.3 in \cite{RaoSongVondra} to $\phi(\lambda):=
(\lambda+m^{2/\alpha})^{\alpha/2}-m$, $\lambda>0$, we can obtain that
the Green kernel $r(x,y):=\int_0^{\infty}p_t(x,y)dt$ of
$X$ satisfies $r(x,y)\asymp (K_{\alpha}(x,y)+K_2(x,y))$, $x,y\in\R^d$,
where $K_{\beta}(x,y):=A(d,\beta)/|x-y|^{d-\beta}$ for $\beta\in]0,2]$. In particular, $X$ is transient and $r(x,x)=\infty$ for $x\in\R^d$.  Note that
$r(x,y)=K_{\alpha}(x,y)$ provided $m=0$.
Let $u$ be a Borel function satisfying $u(x)\Psi(m^{1/\alpha}|x|)/(1+|x|)^{d+\alpha}\in L^1(\R^d)$. For
$\varepsilon>0$ and $x\in\R^d$, we define the modified fractional Laplacian by
$$
\Delta_{\varepsilon}^{\alpha/2,m}u(x):=A(d,-\alpha)\int_{|x-y|>\varepsilon}\frac{u(y)-u(x)}{|x-y|^{d+\alpha}}\Psi(m^{1/\alpha}|x-y|)dy,
$$
and put $\Delta^{\alpha/2,m}u(x):=\lim_{\varepsilon\to0}\Delta_{\varepsilon}^{\alpha/2,m}u(x)$
whenever the limit exists. It is essentially shown in Lemma~3.5 in \cite{BogBya}
(resp.~the remark after Definition~3.7
in \cite{BogBya})
that for any $u\in C_c^2(D)$ (resp.~$u\in C^2(D)$ satisfying
$u(x)\Psi(m^{1/\alpha}|x|)/(1+|x|)^{d+\alpha}\in L^1(\R^d)$),
$\Delta^{\alpha/2,m}u$ always exists in $C(\R^d)$ (resp.~in $C(D)$).
Recall that for $u\in C^2(\R^d)$ with $u(x)\Psi(m^{1/\alpha}|x|)/(1+|x|)^{d+\alpha}\in L^1(\R^d)$, $u$ satisfies \eqref{e:2.1} and \eqref{e:2.2}. Hence, for such $u$ and $\varphi\in C_c^2(D)$,
$\E(u,\varphi)$ is well-defined and the proof of Lemma~2.6 in \cite{CHarm} shows
$$
\int_{\R^d\times\R^d}|u(x)-u(y)||\varphi(x)-\varphi(y)|\frac{\Psi(m^{1/\alpha}|x-y|)dxdy}{|x-y|^{d+\alpha}}<\infty,
$$
which implies $\E(u,\varphi)=(-\Delta^{\alpha/2,m}u,\varphi)$ and
the $\E$-subharmonicity in $D$ of $u$ is equivalent to $\Delta^{\alpha/2,m}u\leq0$ on $D$.

For $\varphi\in C_c(\R^d)$, we set
$$
R^{(\alpha)}\varphi(x):=\int_{\R^d}r(x,y)\varphi(y)dy\quad x\in \R^d.
$$
Then, we see $R^{(\alpha)}\varphi$ is locally bounded on $\R^d$ and
$(R^{(\alpha)}\varphi)(x)
\Psi(m^{1/\alpha}|x|)/(1+|x|)^{d+\alpha}\in L^1(\R^d)$ for such $\varphi$, because of
$r(x,y)\asymp (K_{\alpha}(x,y)+K_2(x,y))$.
Moreover, we see $R^{(\alpha)}\varphi\in\F_{\loc}$ for such $\varphi$. Indeed,
for any relatively compact
 open set $D$ with $\wbar{D}\subset \R^d$,
 $R^{(\alpha)}\varphi$ is
a difference of excessive functions with respect to
 $X^D$ and bounded on $D$, so $R^{(\alpha)}\varphi\in\F_{D,\loc}$
by Theorem~\ref{T:3.9}.
  Since $D$ is arbitrary,
 $R^{(\alpha)}\varphi\in\F_{\loc}$. Thus
$R^{(\alpha)}\varphi$ satisfies
\eqref{e:2.1} and \eqref{e:2.2} for $U,V$ with $\wbar{U}\subset V\subset\wbar{V}\subset \R^d$.
Similarly,
$r(a,\cdot)\in L^{\infty}_{\loc}(\R^d\setminus\{a\})$
satisfies $\int_{\R^d}\frac{r(a,x)
\Psi(m^{1/\alpha}|x|)}{(1+|x|)^{d+\alpha}}dx<\infty$. We can obtain $r(a,\cdot)\in\F_{\R^d\setminus\{a\},\loc}$ in a similar way as above.
Hence $r(a,\cdot)$ satisfies
\eqref{e:2.1} and \eqref{e:2.2} for $U,V$ with $\wbar{U}\subset V\subset\wbar{V}\subset \R^d\setminus\{a\}$.
Note that for $\varphi\in C_c^{\infty}(D)$,
$\Delta^{\alpha/2,m}\varphi=L^{\alpha,m}\varphi$ a.e.~on $\R^d$ and
$R^{(\alpha)}\Delta^{\alpha/2,m}\varphi=-\varphi$ on $\R^d$. Here
$L^{\alpha,m}$ is the $L^2$-generator of $(\E^{\text{\tiny\rm R},\alpha},\F^{\text{\tiny\rm R},\alpha})$.

For $\varphi\in C_c^{\infty}(\R^d\setminus\{a\})$, we then have
\begin{align*}
\E(r(a,\cdot),\varphi)&=-\int_{\R^d}r(a,x)
\Delta^{\alpha/2,m}\varphi(x)dx \\
&=-(R^{(\alpha)}\Delta^{\alpha/2,m}\varphi)(a)=\varphi(a)=0.
\end{align*}
This means the
$\E$-harmonicity  in $\R^d\setminus\{a\}$ of $r(a,\cdot)$.  Similarly, for non-negative
$\psi,\varphi\in C_c^{\infty}(\R^d)$, we have
\begin{align*}
\E(R^{(\alpha)}\psi,\varphi)
=(\psi,-R^{(\alpha)}\Delta^{\alpha/2,m}\varphi)=(\psi,\varphi)\geq0,
\end{align*}
 which implies
the $\E$-superharmonicity of $R^{(\alpha)}\psi$ for
 non-negative $\psi\in C_c^{\infty}(\R^d)$.
\end{example}

\begin{example}[Diffusion process on a locally compact separable metric space]\label{E:4.3}
{\rm Let $(\E,\F)$ be a local regular
Dirichlet form on $L^2(E;m)$, where $E$ is a locally compact separable metric space, and $X$ is its associated
Hunt process. In this case, $X$ has continuous sample paths and so the jumping measure $J$ is null
(cf.~\cite{FOT}). Hence conditions \eqref{e:2.1} and \eqref{e:2.2}  are
automatically satisfied. Let $D$ be an open subset of $E$ and $u$ be a nearly Borel
function on $E$ that is locally bounded in $D$. Then by Theorem~\ref{T:3.11},
$u$ is subharmonic in $D$ if and only if
$u$ is $\E$-subharmonic in $D$.

Now consider the following special case: $E=\R^d$ with $d\geq1$, $m(dx)$ is the
Lebesgue measure $dx$ on $\R^d$, $\F=W^{1,2}(\R^d):=\{u\in L^2(\R^d)\mid\nabla u\in L^2(\R^d)\}$ and
$$
\E(u,v)=\frac{1}{2}\sum_{i,j=1}^d\int_{\R^d}a_{ij}(x)\frac{\partial u(x)}{\partial x_i}
\frac{\partial v(x)}{\partial x_j}dx\quad \text{ for }u,v\in W^{1,2}(\R^d),
$$
where $(a_{i,j}(x))_{1\leq i,j\leq d}$ is a $d\times d$-matrix valued measurable function on
$\R^d$ that is uniformly elliptic and bounded. In literature, $W^{1,2}(\R^d)$ is the Sobolev
space on $\R^d$ of order $(1.2)$. For an open set $D\subset\R^d$, $W^{1,2}(D)$ is
similarly defined as above but with $D$ in lace of $\R^d$. Then $(\E,\F)$ becomes a
regular local Dirichlet form on $L^2(\R^d)$
and its associated Hunt process $X$ is a conservative diffusion on $\R^d$
having jointly
continuous transition density function. Let $D$ be an open set in $\R^d$.
Then by Theorem~\ref{T:3.11}, the following
are equivalent for
 a locally bounded nearly Borel measurable function $u$ on $D$.
\begin{enumerate}
\item\label{item:subharmdiffusion1} $u$ is subharmonic in $D$;
\item\label{item:subharmdiffusion2}
For every relatively compact open subset $U$ of $D$, $u(X_{\tau_U})\in L^1(\bP_x)$ and
$u(x)\leq\bE_x[u(X_{\tau_U})]$ for q.e.~$x\in U$;
$u$ is subharmonic in $D$ in the weak sense;
\item\label{item:subharmdiffusion3} $u\in \F_{D,\loc} = W^{1,2}_{\loc}(D)$ and
$$
\sum_{i,j=1}^d\int_{\R^d}a_{ij}(x)\frac{\partial u(x)}{\partial x_i}
\frac{\partial v(x)}{\partial x_j}dx\leq0\quad\text{ for every }v\in W^{1,2}(D)\cap C_c^+(D).
$$
\end{enumerate}
}
\end{example}

\begin{example}[Diffusions with jumps on $\R^d$]\label{E:4.4}
{\rm
Consider the following Dirichlet form $(\E,\F)$, where $\F =
W^{1,2}(\R^d)$ and for $u,v\in W^{1,2}(\R^d)$
\begin{eqnarray}
\E(u,v):&=&\frac{1}{2}\sum_{i,j=1}^d\int_{\R^d}a_{ij}(x)\frac{\partial u(x)}{\partial x_i}
\frac{\partial v(x)}{\partial x_j}dx \nonumber \\
&&\hskip 0.1truein +\frac{1}{2}\int_{\R^d\times\R^d}(u(x)-u(y))(v(x)-v(y))
   \frac{c(x,y)}{|x-y|^{d+\alpha}}dxdy.  \label{e:4.6}
\end{eqnarray}
Here $d\geq1$ and $(a_{i,j}(x))_{1\leq i,j\leq d}$ is a $d\times d$-matrix valued measurable function on
$\R^d$ that is uniformly elliptic and bounded, $\alpha\in
]0,2[ $ and $c(x, y)$ is
a symmetric function in $(x, y)$ that is bounded between two positive
constants. It is easy to check that $(\E,\F)$ is a regular Dirichlet form on $L^2(\R^d)$.
Its associated symmetric
Hunt process $X$ has both the diffusion and jumping components. Such a process has recently been studied in \cite{CK:apriori}.
Note that when $(a_{i,j}(x))_{1\leq i,j\leq d}$ is the identity matrix and $c(x, y)$ is constant,
the process $X$ is nothing but the
symmetric L\'evy process that is the independent sum of a Brownian motion and a rotationally symmetric
$\alpha$-stable process on $\R^d$.
 It is shown in \cite{CK:apriori} that the Hunt process $X$ associated with the Dirichlet
 form $({\cal E}, W^{1,2}(\R^d))$ given by \eqref{e:4.6}
has strictly positive jointly continuous transition density function
$p_t(x, y)$ and hence is irreducible. Moreover, a sharp two-sided estimate is obtained in
\cite{CK:apriori} for $p_t(x,y)$. In
particular, there is a constant $c > 0$ such that
$$
p_t(x,y)\leq c\left(t^{-d/\alpha}\wedge t^{-d/2} \right)\quad\text{ for any }t>0\text{ and }x,y\in\R^d.
$$
In this example, the jumping measure
$$
J(dxdy)=\frac{c(x,y)}{|x-y|^{d+\alpha}}dxdy.
$$
Hence for any non-empty open set $D\subset \R^d$,
condition \eqref{e:2.1}  is satisfied if and only if $(1 \wedge |x|^{-d-\alpha})u(x)\in L^1(\R^d)$.
 By \cite[Example 2.14]{CHarm}, for this example,
  condition \eqref{e:2.2} is implied by condition \eqref{e:2.1}.
So Theorem~\ref{T:3.11}
 asserts
 that for an open set $D$ and a nearly Borel
  measurable
function $u$ on $\R^d$ that
is locally bounded on $D$ with $(1 \wedge |x|^{-d-\alpha})u(x)\in L^1(\R^d)$,
the following are equivalent
\begin{enumerate}
\item\label{item:subharmdiffusionjump1} $u$ is subharmonic in $D$;
\item\label{item:subharmdiffusionjump2}
For every relatively compact open subset $U$ of $D$, $u(X_{\tau_U})\in L^1(\bP_x)$ and
$u(x)\leq\bE_x[u(X_{\tau_U})]$ for q.e.~$x\in U$;
$u$ is subharmonic in $D$ in the weak sense;
\item\label{item:subharmdiffusionjump3} $u\in \F_{D,\loc} = W^{1,2}_{\loc}(D)$ and
\begin{align*}
\sum_{i,j=1}^d\int_{\R^d}a_{ij}(x)\frac{\partial u(x)}{\partial x_i}
\frac{\partial v(x)}{\partial x_j}dx
+&\int_{\R^d\times\R^d}
(u(x)-u(y))(v(x)-v(y))\frac{c(x,y)}{|x-y|^{d+\alpha}}dxdy
\leq0 \\ &\hspace{3cm}\text{ for every }v\in W^{1,2}(D)\cap C_c^+(D).
\end{align*}
\end{enumerate}

}
\end{example}

\vskip 0.6truein

\noindent {\bf Zhen-Qing Chen}:

Department of Mathematics,
University of Washington, Seattle,
WA 98195, USA.

Email:  {\texttt zchen@math.washington.edu}

\medskip
\noindent {\bf Kazuhiro Kuwae}:

Department of Mathematics and Engineering, Graduate School of
Science and Technology, Kumamoto University, Kumamoto 860-8555, Japan.

 Email: {\texttt kuwae@gpo.kumamoto-u.ac.jp}

\end{document}